\documentclass[11pt]{article}

\usepackage[margin=1in]{geometry}
\usepackage{amsmath,amssymb,mathtools}
\usepackage{newtxtext,newtxmath}
\usepackage{microtype}
\usepackage{setspace}
\usepackage{enumitem}
\providecommand{\tightlist}{\setlength{\itemsep}{0pt}\setlength{\parskip}{0pt}}
\usepackage[hidelinks]{hyperref}
\hypersetup{
  pdftitle={Density-Dependent Operators on Density-Projection Condensation Spaces: Ambient Extensions, Zero-Density Defects, and Stability},
  pdfauthor={Seonghyun Jeon},
  pdfsubject={Functional Analysis},
  pdfkeywords={weighted L2-space, varying Hilbert space, zero-density defect, operator extension}
}

\title{\textbf{Density-Dependent Operators on Density-Projection Condensation Spaces}\\[0.45em]
\large Ambient Extensions, Zero-Density Defects, and Stability}
\author{Seonghyun Jeon}
\date{}

\begin{document}
\maketitle

\begin{abstract}
Let \(X=\mathbb R^n\) be equipped with Lebesgue measure \(\lambda\), and
let

\[
\mathcal R=\left\{\rho\in L^1(X,\lambda):\rho\ge 0\ \lambda\text{-a.e.},\ \int_X\rho\,d\lambda=1\right\}.
\]

Each \(\rho\in\mathcal R\) induces the weighted Hilbert space
\(H_\rho=L^2(X,\rho\,d\lambda)\). Through the alignment isometry

\[
U_\rho[f]_\rho=\sqrt{\rho}\,f,
\]

\(H_\rho\) is identified with the closed observable subspace

\[
V_\rho=\{u\in L^2(X,\lambda):u=0\ \lambda\text{-a.e. on }\{\rho=0\}\}.
\]

The companion density-projection completion theorem identifies the
metric completion of the aligned object space with
\(L^2(X,\lambda)\times\mathcal R\). We study bounded operator families
\(A_\rho:H_\rho\to H_\rho\) whose Hilbert spaces vary with the density.
The aligned operator \(A_\rho^\sharp=U_\rho A_\rho U_\rho^{-1}\) is
defined only on \(V_\rho\). We characterize all bounded extensions of
\(A_\rho^\sharp\) to the ambient space \(L^2(X,\lambda)\) and show that
their collection is an affine space modeled on the bounded operators
from the zero-density defect subspace \(Z_\rho:=V_\rho^\perp\) into the
ambient Hilbert space. The extension that annihilates \(Z_\rho\) is
shown to attain the minimum possible operator norm. We also classify
self-adjoint, positive, and orthogonal-projection extensions and
determine the spectrum of the zero-defect extension.

For a density-dependent extension policy
\(\rho\mapsto\widetilde A_\rho\), we prove that the induced map

\[
(u,\rho)\longmapsto(\widetilde A_\rho u,\rho)
\]

is continuous on the condensation space if and only if the ambient
operator family is strongly continuous in the density. A global
Lipschitz rigidity theorem shows that, for the density-preserving
product dynamics considered here, a linear ambient family can induce a
globally Lipschitz map on the full product space only when it is
independent of \(\rho\). We further characterize strong and
operator-norm convergence of the observable support projections and show
that \(L^1\)-convergence of densities does not, by itself, control
support-dependent operators. Finally, we construct stable
density-dependent multiplication operators and density-weighted
Hilbert-Schmidt integral operators, the latter satisfying a
\(1/2\)-Hölder estimate with respect to the \(L^1\)-distance between
densities.
\end{abstract}

\noindent\textbf{Keywords:} weighted \(L^2\)-space; varying Hilbert space; bounded operator
extension; zero-density defect; strong-operator continuity; support
projection; Hilbert-Schmidt operator.

\section{1. Introduction}\label{introduction}

Let \(X=\mathbb R^n\) with its Borel \(\sigma\)-algebra and Lebesgue
measure \(\lambda\). Every normalized nonnegative density
\(\rho\in L^1(X,\lambda)\) determines the probability measure

\[
d\mu_\rho=\rho\,d\lambda
\]

and the weighted Hilbert space

\[
H_\rho=L^2(X,\mu_\rho).
\]

When \(\rho\) varies, both the Hilbert space and the operators acting on
it may vary. Thus, for

\[
A_\rho:H_\rho\to H_\rho,
\qquad
A_\sigma:H_\sigma\to H_\sigma,
\]

the formal difference \(A_\rho-A_\sigma\) has no intrinsic meaning
unless the two operators are first placed in a common ambient framework.
The same issue occurs at the vector level: elements of \(H_\rho\) and
\(H_\sigma\) cannot be subtracted intrinsically without an
identification of the two fibers.

In the preceding density-projection construction \cite{DPCompletion}, each weighted
space was aligned with the fixed Hilbert space

\[
\mathcal H:=L^2(X,\lambda)
\]

through

\[
U_\rho:H_\rho\to\mathcal H,
\qquad
U_\rho[f]_\rho=\sqrt{\rho}\,f.
\]

This map is unitary from \(H_\rho\) onto

\[
V_\rho:=\{u\in\mathcal H:u=0\ \lambda\text{-a.e. on }\{\rho=0\}\}.
\]

The corresponding density-projection object space carries the metric

\[
d_{DP}\bigl((\rho,\xi),(\sigma,\eta)\bigr)
=
\|U_\rho\xi-U_\sigma\eta\|_2+\|\rho-\sigma\|_1,
\]

and its metric completion was identified isometrically with

\[
\widehat{\mathcal D}^{DP}\cong\mathcal H\times\mathcal R.
\]

The original aligned image consists only of pairs \((u,\rho)\) for which
\(u\) vanishes on the zero-density region of \(\rho\). The completion
also contains pairs whose structural component remains nonzero on such a
region. The orthogonal decomposition

\[
u=u\mathbf 1_{\{\rho>0\}}+u\mathbf 1_{\{\rho=0\}}
\]

separates the observable component from the zero-density defect
component.

This boundary structure creates a new operator-theoretic problem. Given
\(A_\rho\in\mathcal B(H_\rho)\), its aligned realization

\[
A_\rho^\sharp:=U_\rho A_\rho U_\rho^{-1}
\]

acts on \(V_\rho\), but it does not specify an action on

\[
Z_\rho:=V_\rho^\perp.
\]

Since the completion contains the entire ambient space \(\mathcal H\),
an operator on the completion requires an additional choice on
\(Z_\rho\).

The first purpose of this paper is to classify this extension freedom.
We prove that every bounded ambient extension of \(A_\rho^\sharp\) has
the form

\[
\widetilde A_\rho
=
A_\rho^\sharp Q_\rho+B_\rho(I-Q_\rho),
\]

where \(Q_\rho u=\mathbf 1_{\{\rho>0\}}u\) is the orthogonal projection
onto \(V_\rho\), and \(B_\rho:Z_\rho\to\mathcal H\) is an arbitrary
bounded operator. Consequently, the extension space is an affine space
modeled on \(\mathcal B(Z_\rho,\mathcal H)\). In particular, the
intrinsic operator determines the observable action but not the action
on the zero-density defect.

Among all extensions, the choice

\[
\widetilde A_\rho^{\,0}:=A_\rho^\sharp Q_\rho
\]

annihilates the defect subspace. We call it the zero-defect extension.
It is canonical relative to the selected alignment and the requirement
that \(Z_\rho\) be mapped to zero. We prove that it attains the minimum
possible operator norm, although it need not be the unique extension
with that norm.

The second purpose is to study density-dependent extension policies. An
extension policy is a selection \(\rho\mapsto\widetilde A_\rho\)
satisfying \(\widetilde A_\rho|_{V_\rho}=A_\rho^\sharp\). It induces

\[
\mathbb A:\mathcal H\times\mathcal R\to\mathcal H\times\mathcal R,
\qquad
\mathbb A(u,\rho)=(\widetilde A_\rho u,\rho).
\]

We establish an exact continuity criterion: \(\mathbb A\) is continuous
if and only if \(\rho\mapsto\widetilde A_\rho\) is continuous in the
strong operator topology. We also prove a global Lipschitz rigidity
result. If \(\mathbb A\) is globally Lipschitz on the full product
space, then \(\widetilde A_\rho\) must be independent of \(\rho\).

The third purpose is to distinguish support-sensitive and
amplitude-sensitive operator families. The observable projection

\[
Q_\rho=M_{\mathbf 1_{\{\rho>0\}}}
\]

depends only on the positivity set of the density. Even if
\(\rho_k\to\rho\) in \(L^1\), the projections \(Q_{\rho_k}\) need not
converge strongly to \(Q_\rho\). We characterize strong convergence by
local convergence in measure of the positivity sets and show that
operator-norm convergence requires eventual equality of those sets
modulo null sets.

This instability is not universal. We construct multiplication operators
whose coefficients depend continuously on a smoothed density. We also
study the density-weighted integral operators

\[
(\widetilde T_\rho u)(x)
=
\int_X\sqrt{\rho(x)}K(x,y)\sqrt{\rho(y)}u(y)\,d\lambda(y),
\]

where \(K\in L^\infty(X\times X)\). These operators satisfy

\[
\|\widetilde T_\rho-\widetilde T_\sigma\|_{\mathrm{op}}
\le
2\|K\|_\infty\|\rho-\sigma\|_1^{1/2}.
\]

Thus an operator family may remain stable even when the positivity sets
of the densities change discontinuously.

Unlike general convergence theories on varying Hilbert spaces, the
present paper fixes an exact density alignment and studies the extension
data created by the orthogonal defect complement in its explicit metric
completion.

The main contributions are:

\begin{enumerate}
\def\labelenumi{\arabic{enumi}.}
\tightlist
\item
  classification of all bounded ambient extensions of an aligned
  intrinsic operator;
\item
  minimum-norm, structural, and spectral properties of the zero-defect
  extension;
\item
  an exact continuity criterion for induced operators on the
  condensation space;
\item
  a global Lipschitz rigidity theorem;
\item
  strong and operator-norm convergence criteria for observable support
  projections;
\item
  an example showing that two extension policies for the same intrinsic
  identity family can produce different completion dynamics; and
\item
  stable multiplication and Hilbert-Schmidt operator families.
\end{enumerate}

The present paper is restricted to bounded linear operators. Unbounded
operators, density-dependent domains, quadratic forms, semigroup
generators, and resolvent convergence require additional techniques and
are left for future work.

\section{2. Relation to Existing
Frameworks}\label{relation-to-existing-frameworks}

Operators on varying Hilbert spaces have been studied through several
functional-analytic frameworks. Kuwae and Shioya developed a theory for
convergence of Hilbert spaces, bounded operators, quadratic forms,
semigroups, resolvents, and spectral structures \cite{KuwaeShioya2003}. A second
approach introduces identification operators between varying Hilbert
spaces or between those spaces and a common comparison space.
Quantitative compatibility conditions on the identification maps can
then yield resolvent, functional-calculus, semigroup, and spectral
convergence \cite{MugnoloNittkaPost2013,PostZimmer2022}. Varying Hilbert fibers can also be organized
through measurable-field and direct-integral constructions \cite{Yetter2005}.

The present paper does not introduce a general convergence theory for
operators on varying spaces. Instead, it begins with the exact density
alignment

\[
U_\rho:H_\rho\to\mathcal H,
\qquad
U_\rho[f]_\rho=\sqrt{\rho}\,f,
\]

whose range is the concrete subspace \(V_\rho\). The preceding
completion theorem identifies the full completion with
\(\mathcal H\times\mathcal R\) \cite{DPCompletion}. The complement
\(Z_\rho=V_\rho^\perp\) therefore appears as a genuine boundary
component of that completion.

The block decomposition of bounded operators, the strong operator
topology, the Uniform Boundedness Principle, and Hilbert-Schmidt
estimates used below are standard tools \cite{Conway1990,Folland1999,ReedSimon1980}. Likewise, the
existence of bounded extensions from a closed Hilbert subspace is
elementary. The contribution claimed here is narrower: these tools are
applied to the explicit density-aligned completion, and the resulting
extension freedom is interpreted as boundary data carried by the
zero-density defect. More specifically, we identify the following
phenomena:

\begin{enumerate}
\def\labelenumi{\arabic{enumi}.}
\tightlist
\item
  extension freedom is parameterized precisely by operators acting on
  the zero-density defect subspace;
\item
  the same intrinsic operator family may generate different dynamics
  under different defect-extension policies;
\item
  support-sensitive policies may be unstable under \(L^1\)-convergence
  of densities; and
\item
  amplitude-weighted families can remain stable across support
  degeneration.
\end{enumerate}

Accordingly, the results should be viewed as an analysis of bounded
operator extension and stability within a particular density-dependent
completion, rather than as a replacement for general theories of varying
Hilbert spaces.

\section{3. Preliminaries and
Notation}\label{preliminaries-and-notation}

Throughout the paper, \(X=\mathbb R^n\) is equipped with its Borel
\(\sigma\)-algebra and Lebesgue measure \(\lambda\). Functions are
complex-valued unless otherwise stated. Inner products are linear in the
first variable and conjugate-linear in the second. We write

\[
\mathcal H=L^2(X,\lambda)
\]

and denote its norm by \(\|\cdot\|_2\). The operator norm is denoted by
\(\|\cdot\|_{\mathrm{op}}\).

\subsection{3.1. Normalized densities}\label{normalized-densities}

Define

\[
\mathcal R
=
\left\{\rho\in L^1(X,\lambda):\rho\ge0\ \lambda\text{-a.e.},\ \int_X\rho\,d\lambda=1\right\}.
\]

The space \(\mathcal R\) is equipped with the \(L^1\)-metric. It is
closed in \(L^1(X,\lambda)\) and therefore complete. For
\(\rho\in\mathcal R\), define

\[
\mu_\rho(A)=\int_A\rho\,d\lambda.
\]

Then \(\mu_\rho\) is a probability measure and
\(d\mu_\rho/d\lambda=\rho\).

\subsection{3.2. Weighted Hilbert spaces and
alignment}\label{weighted-hilbert-spaces-and-alignment}

For \(\rho\in\mathcal R\), let

\[
H_\rho=L^2(X,\mu_\rho).
\]

An element is denoted by \([f]_\rho\), with norm and inner product

\[
\|[f]_\rho\|_{H_\rho}^2=\int_X|f|^2\rho\,d\lambda,
\]

\[
\langle[f]_\rho,[g]_\rho\rangle_{H_\rho}=\int_Xf\overline g\,\rho\,d\lambda.
\]

Define

\[
U_\rho:H_\rho\to\mathcal H,
\qquad
U_\rho[f]_\rho=\sqrt{\rho}\,f.
\]

This is a well-defined linear isometry with range

\[
V_\rho
=
\{u\in\mathcal H:u=0\ \lambda\text{-a.e. on }\{\rho=0\}\}.
\]

Thus \(U_\rho:H_\rho\to V_\rho\) is unitary.

Define

\[
Z_\rho:=V_\rho^\perp
=
\{u\in\mathcal H:u=0\ \lambda\text{-a.e. on }\{\rho>0\}\}.
\]

Then

\[
\mathcal H=V_\rho\oplus Z_\rho.
\]

The associated orthogonal projections are

\[
Q_\rho u=\mathbf1_{\{\rho>0\}}u,
\qquad
(I-Q_\rho)u=\mathbf1_{\{\rho=0\}}u.
\]

The positivity and zero sets depend on the representative only by null
sets, so \(V_\rho\), \(Z_\rho\), and \(Q_\rho\) are well defined for the
\(L^1\)-class of \(\rho\).

\subsection{3.3. The condensation space}\label{the-condensation-space}

The intrinsic object space is

\[
\mathcal D^{DP}
=
\bigsqcup_{\rho\in\mathcal R}(\{\rho\}\times H_\rho)
\]

with metric

\[
d_{DP}\bigl((\rho,\xi),(\sigma,\eta)\bigr)
=
\|U_\rho\xi-U_\sigma\eta\|_2+\|\rho-\sigma\|_1.
\]

The alignment

\[
J(\rho,\xi)=(U_\rho\xi,\rho)
\]

is isometric. By the image and completion theorems of the companion
manuscript \cite[Theorems~4.1 and~4.3]{DPCompletion}, its image is

\[
J(\mathcal D^{DP})
=
\{(u,\rho)\in\mathcal H\times\mathcal R:u\in V_\rho\}.
\]

Its metric completion is

\[
\widehat{\mathcal D}^{DP}\cong\mathcal H\times\mathcal R
\]

with

\[
d_\oplus((u,\rho),(v,\sigma))
=
\|u-v\|_2+\|\rho-\sigma\|_1.
\]

For a completed point \((u,\rho)\), define

\[
u_{\mathrm{obs}}:=Q_\rho u,
\qquad
u_{\mathrm{def}}:=(I-Q_\rho)u.
\]

Then \(u=u_{\mathrm{obs}}+u_{\mathrm{def}}\) orthogonally. By {[}5,
Proposition 5.1{]}, \((u,\rho)\) lies in the original aligned image if
and only if \(u_{\mathrm{def}}=0\).

\section{4. Density-Dependent
Operators}\label{density-dependent-operators}

\textbf{Definition 4.1 (Density-dependent operator family).} A bounded
density-dependent operator family is a collection

\[
\mathcal A=\{A_\rho\}_{\rho\in\mathcal R},
\qquad
A_\rho\in\mathcal B(H_\rho).
\]

\textbf{Definition 4.2 (Aligned operator).} For each
\(\rho\in\mathcal R\), define

\[
A_\rho^\sharp:=U_\rho A_\rho U_\rho^{-1}:V_\rho\to V_\rho.
\]

Because \(U_\rho\) is unitary,

\[
\|A_\rho^\sharp\|_{\mathcal B(V_\rho)}
=
\|A_\rho\|_{\mathcal B(H_\rho)}.
\]

Unitary equivalence preserves self-adjointness, positivity, normality,
idempotence, compactness, and spectrum. These statements concern the
operator on the fiber and its aligned realization on \(V_\rho\); they do
not yet specify any action on \(Z_\rho\).

\textbf{Definition 4.3 (Fiber-compatible ambient family).} A family

\[
\{\widetilde A_\rho\}_{\rho\in\mathcal R}\subset\mathcal B(\mathcal H)
\]

is fiber-compatible if

\[
\widetilde A_\rho(V_\rho)\subseteq V_\rho
\]

for every \(\rho\in\mathcal R\).

Every fiber-compatible family determines an intrinsic family by

\[
A_\rho
=
U_\rho^{-1}(\widetilde A_\rho|_{V_\rho})U_\rho.
\]

Conversely, an intrinsic family determines an aligned operator only on
\(V_\rho\). Extending it to \(\mathcal H\) is the central problem of the
next section.

\section{5. Ambient Extension Spaces}\label{ambient-extension-spaces}

\textbf{Definition 5.1 (Ambient extension).} An operator
\(\widetilde A\in\mathcal B(\mathcal H)\) is an ambient extension of
\(A_\rho\) if

\[
\widetilde A|_{V_\rho}=A_\rho^\sharp.
\]

The collection of all such extensions is denoted by
\(\operatorname{Ext}_\rho(A_\rho)\).

Define the zero-defect extension

\[
\widetilde A_\rho^{\,0}:=A_\rho^\sharp Q_\rho.
\]

\textbf{Theorem 5.2 (Affine structure of the extension space).} For each
\(\rho\in\mathcal R\),

\[
\operatorname{Ext}_\rho(A_\rho)
=
\widetilde A_\rho^{\,0}+\mathcal N_\rho,
\]

where

\[
\mathcal N_\rho
=
\{N\in\mathcal B(\mathcal H):NQ_\rho=0\}.
\]

Moreover,

\[
\mathcal N_\rho\cong\mathcal B(Z_\rho,\mathcal H)
\]

isometrically. Hence \(\operatorname{Ext}_\rho(A_\rho)\) is an affine
space modeled on \(\mathcal B(Z_\rho,\mathcal H)\).

\textbf{Proof.} Let \(\widetilde A\in\operatorname{Ext}_\rho(A_\rho)\).
For \(u\in\mathcal H\),

\[
\widetilde Au
=
\widetilde A Q_\rho u+\widetilde A(I-Q_\rho)u
=
A_\rho^\sharp Q_\rho u+\widetilde A(I-Q_\rho)u.
\]

Thus

\[
\widetilde A=\widetilde A_\rho^{\,0}+N,
\qquad
N:=\widetilde A(I-Q_\rho),
\]

and \(NQ_\rho=0\). Conversely, if \(NQ_\rho=0\), then for
\(v\in V_\rho\),

\[
(\widetilde A_\rho^{\,0}+N)v=A_\rho^\sharp v,
\]

so \(\widetilde A_\rho^{\,0}+N\) is an extension.

For \(B\in\mathcal B(Z_\rho,\mathcal H)\), the operator \(B(I-Q_\rho)\)
belongs to \(\mathcal N_\rho\). Conversely, each \(N\in\mathcal N_\rho\)
satisfies \(N=N(I-Q_\rho)\) and is determined by \(N|_{Z_\rho}\).
Finally,

\[
\|B(I-Q_\rho)\|_{\mathrm{op}}=\|B\|_{\mathrm{op}},
\]

because \(I-Q_\rho\) is the orthogonal projection onto \(Z_\rho\).
\(\square\)

\textbf{Corollary 5.3 (Extension uniqueness).} The ambient extension is
unique if and only if \(Z_\rho=\{0\}\). Equivalently, uniqueness holds
if and only if \(\rho>0\) almost everywhere.

\textbf{Proof.} By Theorem 5.2, the extension freedom is modeled on
\(\mathcal B(Z_\rho,\mathcal H)\), which is trivial exactly when
\(Z_\rho=\{0\}\). \(\square\)

\textbf{Corollary 5.4 (Uniqueness under defect annihilation).} The
zero-defect extension \(\widetilde A_\rho^{\,0}\) is the unique ambient
extension satisfying

\[
\widetilde A z=0
\qquad
\text{for every }z\in Z_\rho.
\]

\textbf{Proof.} In the representation of Theorem 5.2, the restriction to
\(Z_\rho\) is the free operator \(B\). The stated condition forces
\(B=0\). \(\square\)

Thus extension nonuniqueness is determined exactly by the zero-density
defect space, while the additional requirement of defect annihilation
selects a unique extension.

\section{6. Canonical and Structure-Preserving
Extensions}\label{canonical-and-structure-preserving-extensions}

\subsection{6.1. Minimum operator norm}\label{minimum-operator-norm}

\textbf{Theorem 6.1 (Minimum-norm zero-defect extension).} For every
\(\widetilde A\in\operatorname{Ext}_\rho(A_\rho)\),

\[
\|\widetilde A\|_{\mathrm{op}}
\ge
\|A_\rho\|_{\mathcal B(H_\rho)}.
\]

The zero-defect extension satisfies

\[
\|\widetilde A_\rho^{\,0}\|_{\mathrm{op}}
=
\|A_\rho\|_{\mathcal B(H_\rho)}.
\]

Hence it attains the minimum possible operator norm among all bounded
ambient extensions.

\textbf{Proof.} Every extension agrees with \(A_\rho^\sharp\) on
\(V_\rho\), so

\[
\|\widetilde A\|_{\mathrm{op}}
\ge
\|\widetilde A|_{V_\rho}\|_{
\mathcal B(V_\rho)}
=
\|A_\rho^\sharp\|_{
\mathcal B(V_\rho)}
=
\|A_\rho\|_{
\mathcal B(H_\rho)}.
\]

On the other hand,

\[
\widetilde A_\rho^{\,0}=A_\rho^\sharp Q_\rho
\]

and \(Q_\rho\) is a contraction. Therefore

\[
\|\widetilde A_\rho^{\,0}\|_{\mathrm{op}}
\le
\|A_\rho^\sharp\|_{\mathrm{op}}.
\]

Restriction to \(V_\rho\) gives the reverse inequality. \(\square\)

The zero-defect extension need not be the unique extension attaining the
minimum norm. The word \emph{canonical} is used only relative to the
chosen alignment and the requirement that \(Z_\rho\) be annihilated.

\subsection{6.2. Block representation and preserved
structure}\label{block-representation-and-preserved-structure}

Relative to \(\mathcal H=V_\rho\oplus Z_\rho\), every extension has the
block form

\[
\widetilde A
=
\begin{pmatrix}
A_\rho^\sharp & C\\
0 & D
\end{pmatrix},
\]

where \(C\in\mathcal B(Z_\rho,V_\rho)\) and \(D\in\mathcal B(Z_\rho)\).

\textbf{Theorem 6.2 (Structure-preserving extensions).}

\begin{enumerate}
\def\labelenumi{\arabic{enumi}.}
\tightlist
\item
  \(\widetilde A\) is self-adjoint if and only if \[
  A_\rho^\sharp=(A_\rho^\sharp)^*,\qquad C=0,\qquad D=D^*.
  \]
\item
  \(\widetilde A\) is positive if and only if \[
  A_\rho^\sharp\ge0,\qquad C=0,\qquad D\ge0.
  \]
\item
  \(\widetilde A\) is an orthogonal projection if and only if \(C=0\)
  and both \(A_\rho^\sharp\) and \(D\) are orthogonal projections on
  \(V_\rho\) and \(Z_\rho\), respectively.
\end{enumerate}

\textbf{Proof.} The adjoint block matrix is

\[
\widetilde A^*
=
\begin{pmatrix}
(A_\rho^\sharp)^* & 0\\
C^* & D^*
\end{pmatrix}.
\]

Equality \(\widetilde A=\widetilde A^*\) yields the first statement. If
\(\widetilde A\) is positive, then it is self-adjoint. Comparing the
upper-right and lower-left blocks of \(\widetilde A\) and
\(\widetilde A^*\) therefore forces \(C=0\) and \(D=D^*\). With \(C=0\),
for every \(v\in V_\rho\) and \(z\in Z_\rho\) one has

\[
\langle \widetilde A(v+z),v+z\rangle
=
\langle A_\rho^\sharp v,v\rangle
+
\langle Dz,z\rangle.
\]

Thus \(\widetilde A\) is positive exactly when both diagonal blocks are
positive. An orthogonal projection is a self-adjoint idempotent, so the
third statement follows by applying the first statement and checking the
two diagonal blocks. \(\square\)

\subsection{6.3. Spectrum of the zero-defect
extension}\label{spectrum-of-the-zero-defect-extension}

The zero-defect extension has the direct-sum form

\[
\widetilde A_\rho^{\,0}=A_\rho^\sharp\oplus0_{Z_\rho}.
\]

\textbf{Proposition 6.3.}

\[
\sigma(\widetilde A_\rho^{\,0})
=
\begin{cases}
\sigma(A_\rho),&Z_\rho=\{0\},\\
\sigma(A_\rho)\cup\{0\},&Z_\rho\ne\{0\}.
\end{cases}
\]

\textbf{Proof.} \(A_\rho\) and \(A_\rho^\sharp\) are unitarily
equivalent. The spectrum of a direct sum is the union of the spectra of
its summands, and \(\sigma(0_{Z_\rho})=\{0\}\) when \(Z_\rho\ne\{0\}\).
\(\square\)

\textbf{Proposition 6.4 (Kernel contribution of the defect).} The kernel
of the zero-defect extension is

\[
\ker(\widetilde A_\rho^{\,0})
=
U_\rho(\ker A_\rho)\oplus Z_\rho.
\]

In particular, if \(Z_\rho\ne\{0\}\), then \(0\) belongs to the point
spectrum of \(\widetilde A_\rho^{\,0}\). If \(A_\rho\) is injective, the
entire kernel of the extension is exactly \(Z_\rho\).

\textbf{Proof.} Relative to \(\mathcal H=V_\rho\oplus Z_\rho\),

\[
\widetilde A_\rho^{\,0}=A_\rho^\sharp\oplus0_{Z_\rho}.
\]

Hence

\[
\ker(\widetilde A_\rho^{\,0})
=
\ker(A_\rho^\sharp)\oplus Z_\rho.
\]

Unitary equivalence gives \(\ker(A_\rho^\sharp)=U_\rho(\ker A_\rho)\).
\(\square\)

If \(A_\rho\) is compact, self-adjoint, positive, normal, or idempotent,
then the zero-defect extension has the corresponding property.

\section{7. Extension Policies and
Continuity}\label{extension-policies-and-continuity}

\textbf{Definition 7.1 (Extension policy).} An extension policy is a
selection

\[
\rho\longmapsto\widetilde A_\rho\in\operatorname{Ext}_\rho(A_\rho)
\]

for every \(\rho\in\mathcal R\).

It induces

\[
\mathbb A:\mathcal H\times\mathcal R\to\mathcal H\times\mathcal R,
\qquad
\mathbb A(u,\rho)=(\widetilde A_\rho u,\rho).
\]

\textbf{Definition 7.2 (Strong continuity in the density).} The policy
is strongly continuous if \(\rho_k\to\rho\) in \(L^1\) implies

\[
\widetilde A_{\rho_k}u\to\widetilde A_\rho u
\]

in \(\mathcal H\) for every fixed \(u\in\mathcal H\).

\textbf{Theorem 7.3 (Exact continuity criterion).} The induced map
\(\mathbb A\) is continuous on \(\mathcal H\times\mathcal R\) if and
only if \(\rho\mapsto\widetilde A_\rho\) is strongly continuous.

\textbf{Proof.} If \(\mathbb A\) is continuous and \(\rho_k\to\rho\),
then \((u,\rho_k)\to(u,\rho)\) for every fixed \(u\), so
\(\widetilde A_{\rho_k}u\to\widetilde A_\rho u\).

Conversely, suppose the policy is strongly continuous and
\((u_k,\rho_k)\to(u,\rho)\). For every fixed \(x\in\mathcal H\), the
sequence \(\widetilde A_{\rho_k}x\) converges, and hence is bounded. By
the Uniform Boundedness Principle \cite{Conway1990},

\[
\sup_k\|\widetilde A_{\rho_k}\|_{\mathrm{op}}<\infty.
\]

Therefore

\[
\begin{aligned}
\|\widetilde A_{\rho_k}u_k-\widetilde A_\rho u\|_2
&\le
\|\widetilde A_{\rho_k}(u_k-u)\|_2\\
&\quad+
\|(\widetilde A_{\rho_k}-\widetilde A_\rho)u\|_2.
\end{aligned}
\]

The first term tends to zero by uniform boundedness and \(u_k\to u\),
while the second tends to zero by strong continuity. \(\square\)

\textbf{Proposition 7.4 (Local Lipschitz criterion).} Suppose

\[
\sup_{\rho\in\mathcal R}\|\widetilde A_\rho\|_{\mathrm{op}}\le M
\]

and

\[
\|\widetilde A_\rho-\widetilde A_\sigma\|_{\mathrm{op}}
\le
L\|\rho-\sigma\|_1.
\]

Then \(\mathbb A\) is Lipschitz on each set
\(\{(u,\rho):\|u\|_2\le R\}\). More precisely,

\[
\begin{aligned}
d_\oplus(\mathbb A(u,\rho),\mathbb A(v,\sigma))
&\le
M\|u-v\|_2\\
&\quad+(LR+1)\|\rho-\sigma\|_1.
\end{aligned}
\]

\textbf{Proof.} Use

\[
\|\widetilde A_\rho u-\widetilde A_\sigma v\|_2
\le
\|\widetilde A_\rho(u-v)\|_2
+
\|(\widetilde A_\rho-\widetilde A_\sigma)v\|_2.
\]

Then add the density distance. \(\square\)

\textbf{Theorem 7.5 (Global Lipschitz rigidity).} The map \(\mathbb A\)
is globally Lipschitz on \(\mathcal H\times\mathcal R\) if and only if
there exists a single \(T\in\mathcal B(\mathcal H)\) such that
\(\widetilde A_\rho=T\) for every \(\rho\in\mathcal R\).

\textbf{Proof.} Suppose \(\mathbb A\) is globally \(C\)-Lipschitz.
Because \(\mathcal R\) contains distinct densities and the density
coordinate is preserved, necessarily \(C\ge1\). Fix \(\rho,\sigma,u\)
and \(t>0\). Comparing \((tu,\rho)\) and \((tu,\sigma)\) gives

\[
t\|(\widetilde A_\rho-\widetilde A_\sigma)u\|_2
+
\|\rho-\sigma\|_1
\le
C\|\rho-\sigma\|_1.
\]

Hence

\[
t\|(\widetilde A_\rho-\widetilde A_\sigma)u\|_2
\le
(C-1)\|\rho-\sigma\|_1.
\]

Dividing by \(t>0\) gives

\[
\|(\widetilde A_\rho-\widetilde A_\sigma)u\|_2
\le
\frac{C-1}{t}\|\rho-\sigma\|_1.
\]

The right-hand side tends to zero as \(t\to\infty\). Therefore
\((\widetilde A_\rho-\widetilde A_\sigma)u=0\). Since \(u\) is
arbitrary, \(\widetilde A_\rho=\widetilde A_\sigma\).

Conversely, if \(\widetilde A_\rho=T\) for all \(\rho\), then

\[
d_\oplus(\mathbb A(u,\rho),\mathbb A(v,\sigma))
\le
\max\{\|T\|_{\mathrm{op}},1\}
 d_\oplus((u,\rho),(v,\sigma)).
\]

\(\square\)

\section{8. Support Projections and Density
Degeneration}\label{support-projections-and-density-degeneration}

Write \(E_\rho=\{\rho>0\}\). Then \(Q_\rho=M_{\mathbf1_{E_\rho}}\).

\textbf{Theorem 8.1 (Strong convergence of support projections).} For
\(\rho_k,\rho\in\mathcal R\), the following are equivalent:

\begin{enumerate}
\def\labelenumi{\arabic{enumi}.}
\tightlist
\item
  \(Q_{\rho_k}\to Q_\rho\) strongly on \(\mathcal H\);
\item
  for every measurable \(B\subset X\) with \(\lambda(B)<\infty\), \[
  \lambda((E_{\rho_k}\triangle E_\rho)\cap B)\to0;
  \]
\item
  \(\mathbf1_{E_{\rho_k}}\to\mathbf1_{E_\rho}\) locally in measure.
\end{enumerate}

\textbf{Proof.} If strong convergence holds, apply \(Q_{\rho_k}-Q_\rho\)
to \(\mathbf1_B\):

\[
\|(Q_{\rho_k}-Q_\rho)\mathbf1_B\|_2^2
=
\lambda((E_{\rho_k}\triangle E_\rho)\cap B).
\]

This proves (1)\(\Rightarrow\)(2). For indicator functions, (2) and (3)
are equivalent.

Assume (2). Let \(u\in L^2(X,\lambda)\) and choose a bounded function
\(v\) supported on a finite-measure set \(B\) such that
\(\|u-v\|_2<\varepsilon\). Since \(Q_{\rho_k}-Q_\rho\) is multiplication
by a function with absolute value at most one,

\[
\begin{aligned}
\|(Q_{\rho_k}-Q_\rho)u\|_2
&\le
\|u-v\|_2
+
\|(Q_{\rho_k}-Q_\rho)v\|_2\\
&\le
\varepsilon+
\|v\|_\infty
\lambda((E_{\rho_k}\triangle E_\rho)\cap B)^{1/2}.
\end{aligned}
\]

The second term tends to zero, and \(\varepsilon\) is arbitrary.
\(\square\)

\textbf{Theorem 8.2 (Operator-norm convergence).}

\[
\|Q_{\rho_k}-Q_\rho\|_{\mathrm{op}}
=
\|\mathbf1_{E_{\rho_k}}-\mathbf1_{E_\rho}\|_\infty
\in\{0,1\}.
\]

Consequently, \(Q_{\rho_k}\to Q_\rho\) in operator norm if and only if
\(E_{\rho_k}=E_\rho\) almost everywhere for all sufficiently large
\(k\).

\textbf{Proof.} The difference is a multiplication operator whose
multiplier takes values in \(\{-1,0,1\}\). Its operator norm is
therefore either zero or one. \(\square\)

\textbf{Example 8.3 (Support degeneration).} Suppose
\(\lambda(\{\rho=0\})>0\), choose \(\gamma\in\mathcal R\) with
\(\gamma>0\) almost everywhere, and define, for \(k\ge2\),

\[
\rho_k=\left(1-\frac1k\right)\rho+\frac1k\gamma.
\]

Then \(\rho_k\to\rho\) in \(L^1\), but \(Q_{\rho_k}=I\) for every \(k\),
whereas \(Q_\rho\ne I\). Thus \(L^1\)-convergence alone does not imply
strong or operator-norm convergence of support projections.

\textbf{Corollary 8.4 (Strong continuity points of the support policy).}
The map

\[
\mathcal R\longrightarrow\mathcal B(\mathcal H),
\qquad
\rho\longmapsto Q_\rho,
\]

is continuous at \(\rho\) in the strong operator topology if and only if

\[
\rho>0
\qquad
\lambda\text{-a.e.}
\]

\textbf{Proof.} If \(\lambda(\{\rho=0\})>0\), Example 8.3 provides a
sequence \(\rho_k\to\rho\) for which \(Q_{\rho_k}=I\) but
\(Q_\rho\ne I\), so strong continuity fails.

Conversely, suppose \(\rho>0\) almost everywhere and \(\rho_k\to\rho\)
in \(L^1\). Let \(B\) have finite measure and set

\[
F_k=B\cap\{\rho_k=0\}.
\]

Because \(\rho_k=0\) on \(F_k\),

\[
\int_{F_k}\rho\,d\lambda
\le
\|\rho_k-\rho\|_1.
\]

For \(\delta>0\),

\[
\lambda(F_k)
\le
\lambda(B\cap\{0<\rho\le\delta\})
+
\delta^{-1}\|\rho_k-\rho\|_1.
\]

Since \(\rho>0\) almost everywhere on \(B\), the first term tends to
zero as \(\delta\downarrow0\). Hence \(\lambda(F_k)\to0\). Because
\(E_\rho=X\) modulo null sets, Theorem 8.1 gives
\(Q_{\rho_k}\to Q_\rho=I\) strongly. \(\square\)

\textbf{Corollary 8.5 (Nowhere norm continuity).} The map
\(\rho\mapsto Q_\rho\) is nowhere continuous in the operator norm.

\textbf{Proof.} If \(\lambda(\{\rho=0\})>0\), the sequence in Example
8.3 satisfies

\[
\|Q_{\rho_k}-Q_\rho\|_{\mathrm{op}}=1.
\]

Now suppose \(\rho>0\) almost everywhere. The probability measure
\(\mu_\rho=\rho\lambda\) is non-atomic. Choose measurable sets \(F_k\)
such that

\[
0<m_k:=\mu_\rho(F_k)<\frac1k.
\]

Define

\[
\sigma_k
=
\frac{\rho\mathbf1_{F_k^c}}{1-m_k}.
\]

Then \(\sigma_k\in\mathcal R\) and

\[
\|\sigma_k-\rho\|_1=2m_k\longrightarrow0.
\]

However, \(Q_{\sigma_k}\) vanishes on \(F_k\) while \(Q_\rho=I\) modulo
null sets. Therefore

\[
\|Q_{\sigma_k}-Q_\rho\|_{\mathrm{op}}=1
\]

for every \(k\). \(\square\)

\section{9. Intrinsic Underdetermination and Extension
Policies}\label{intrinsic-underdetermination-and-extension-policies}

\textbf{Example 9.1 (Two policies for the identity family).} Let

\[
A_\rho=I_{H_\rho}
\]

for every \(\rho\). Then \(A_\rho^\sharp=I_{V_\rho}\).

The zero-defect policy gives

\[
\widetilde A_\rho^{\,0}=Q_\rho,
\]

and, by Corollary 8.4, the induced map
\((u,\rho)\mapsto(Q_\rho u,\rho)\) is discontinuous at every density
whose zero set has positive measure.

A second valid policy is

\[
\widetilde A_\rho^{\,1}=I_{\mathcal H}.
\]

Indeed, \(I_{\mathcal H}|_{V_\rho}=I_{V_\rho}\). The induced map is the
identity on \(\mathcal H\times\mathcal R\) and is globally
\(1\)-Lipschitz.

\textbf{Proposition 9.2.} The intrinsic identity family admits both an
extension policy that is discontinuous at every density with a
nontrivial zero region and a globally Lipschitz extension policy. Hence
the intrinsic operators alone do not determine the dynamics of the
zero-density defect.

The example separates the intrinsic observable action from the selected
defect action. The first is determined by \(A_\rho\); the second must be
supplied by the extension policy.

\section{10. Stable Multiplication
Families}\label{stable-multiplication-families}

Let \(a_\rho\in L^\infty(X,\lambda)\) and define

\[
\widetilde M_\rho u=a_\rho u.
\]

Multiplication preserves \(V_\rho\), so the family is fiber-compatible.
The corresponding intrinsic operator is

\[
M_\rho[f]_\rho=[a_\rho f]_\rho.
\]

\textbf{Proposition 10.1.} Suppose

\[
\|a_\rho\|_\infty\le M
\]

and

\[
\|a_\rho-a_\sigma\|_\infty
\le
L\|\rho-\sigma\|_1.
\]

Then

\[
\mathbb M(u,\rho)=(a_\rho u,\rho)
\]

is locally Lipschitz on the condensation space.

\textbf{Proof.} Multiplication operators satisfy

\[
\|\widetilde M_\rho\|_{\mathrm{op}}=\|a_\rho\|_\infty,
\qquad
\|\widetilde M_\rho-\widetilde M_\sigma\|_{\mathrm{op}}
=
\|a_\rho-a_\sigma\|_\infty.
\]

Apply Proposition 7.4. \(\square\)

\subsection{10.1. Smoothed density
multipliers}\label{smoothed-density-multipliers}

Let \(\kappa\in L^\infty(\mathbb R^n)\) and let
\(\phi:\mathbb C\to\mathbb C\) be bounded and Lipschitz. Define

\[
(\kappa*\rho)(x)
:=
\int_X\kappa(x-y)\rho(y)\,d\lambda(y),
\qquad
 a_\rho(x)=\phi((\kappa*\rho)(x)).
\]

Then

\[
\|\kappa*(\rho-\sigma)\|_\infty
\le
\|\kappa\|_\infty\|\rho-\sigma\|_1,
\]

so

\[
\|a_\rho-a_\sigma\|_\infty
\le
\operatorname{Lip}(\phi)\|\kappa\|_\infty\|\rho-\sigma\|_1.
\]

Thus the associated multiplication family is operator-norm Lipschitz in
the density.

\section{11. Density-Weighted Hilbert-Schmidt
Operators}\label{density-weighted-hilbert-schmidt-operators}

Let \(K\in L^\infty(X\times X)\) be represented by a jointly measurable
kernel, and define

\[
k_\rho(x,y)=\sqrt{\rho(x)}K(x,y)\sqrt{\rho(y)}.
\]

Since

\[
\begin{aligned}
\|k_\rho\|_{L^2(X\times X)}^2
&=
\int_X\int_X\rho(x)|K(x,y)|^2\rho(y)\,d\lambda(y)d\lambda(x)\\
&\le
\|K\|_\infty^2,
\end{aligned}
\]

\(k_\rho\) defines a Hilbert-Schmidt operator

\[
(\widetilde T_\rho u)(x)
=
\int_Xk_\rho(x,y)u(y)\,d\lambda(y).
\]

Because of the factor \(\sqrt{\rho(x)}\), one has
\(\widetilde T_\rho(\mathcal H)\subseteq V_\rho\), so the family is
fiber-compatible.

The corresponding intrinsic operator is

\[
T_\rho[f]_\rho
=
\left[x\mapsto\int_XK(x,y)f(y)\,d\mu_\rho(y)\right]_\rho.
\]

\textbf{Proposition 11.1.} The operator \(T_\rho\) is bounded,

\[
\|T_\rho\|_{\mathcal B(H_\rho)}\le\|K\|_\infty,
\]

and

\[
U_\rho T_\rho=\widetilde T_\rho U_\rho.
\]

\textbf{Proof.} The definition is independent of the representative of
\([f]_\rho\), and measurability follows from the standard
parameter-integral theorem. Since \(\mu_\rho\) is a probability measure,

\[
\left|\int_XK(x,y)f(y)\,d\mu_\rho(y)\right|
\le
\|K\|_\infty\|f\|_{L^1(\mu_\rho)}
\le
\|K\|_\infty\|f\|_{H_\rho}.
\]

The output is bounded and hence belongs to \(H_\rho\). The intertwining
identity follows by direct calculation. \(\square\)

\textbf{Theorem 11.2 (Square-root-density stability).} For all
\(\rho,\sigma\in\mathcal R\),

\[
\|\widetilde T_\rho-\widetilde T_\sigma\|_{\mathrm{op}}
\le
2\|K\|_\infty\|\sqrt\rho-\sqrt\sigma\|_2.
\]

Consequently,

\[
\|\widetilde T_\rho-\widetilde T_\sigma\|_{\mathrm{op}}
\le
2\|K\|_\infty\|\rho-\sigma\|_1^{1/2}.
\]

\textbf{Proof.} Decompose

\[
\begin{aligned}
k_\rho(x,y)-k_\sigma(x,y)
&=
(\sqrt{\rho(x)}-\sqrt{\sigma(x)})K(x,y)\sqrt{\rho(y)}\\
&\quad+
\sqrt{\sigma(x)}K(x,y)(\sqrt{\rho(y)}-\sqrt{\sigma(y)}).
\end{aligned}
\]

Each term has \(L^2(X\times X)\)-norm at most

\[
\|K\|_\infty\|\sqrt\rho-\sqrt\sigma\|_2,
\]

because \(\|\sqrt\rho\|_2=\|\sqrt\sigma\|_2=1\). Hence

\[
\|k_\rho-k_\sigma\|_2
\le
2\|K\|_\infty\|\sqrt\rho-\sqrt\sigma\|_2.
\]

The operator norm is bounded by the Hilbert-Schmidt norm. Finally,

\[
|\sqrt a-\sqrt b|^2\le|a-b|,
\qquad a,b\ge0,
\]

which gives

\[
\|\sqrt\rho-\sqrt\sigma\|_2^2\le\|\rho-\sigma\|_1.
\]

\(\square\)

\textbf{Remark 11.3 (\(1/2\)-Hölder continuity).} The second estimate in
Theorem 11.2 states precisely that

\[
\rho\longmapsto\widetilde T_\rho
\]

is \(1/2\)-Hölder continuous from \((\mathcal R,\|\cdot\|_1)\) into
\(\mathcal B(\mathcal H)\) equipped with the operator norm. This
contrasts sharply with the support policy \(\rho\mapsto Q_\rho\), which
is nowhere continuous in operator norm by Corollary 8.5. Thus support
degeneration may destroy continuity for support-sensitive policies while
remaining compatible with quantitative operator-norm stability for
amplitude-weighted families.

\textbf{Corollary 11.4.} The map

\[
\mathbb T(u,\rho)=(\widetilde T_\rho u,\rho)
\]

is continuous. On \(\{(u,\rho):\|u\|_2\le R\}\),

\[
\begin{aligned}
d_\oplus(\mathbb T(u,\rho),\mathbb T(v,\sigma))
&\le
\|K\|_\infty\|u-v\|_2\\
&\quad+
2R\|K\|_\infty\|\rho-\sigma\|_1^{1/2}\\
&\quad+
\|\rho-\sigma\|_1.
\end{aligned}
\]

This family remains operator-norm continuous even when positivity sets
degenerate. Unlike \(Q_\rho\), it depends on the density amplitude
through the factors \(\sqrt{\rho(x)}\sqrt{\rho(y)}\), which vanish
continuously as the density degenerates.

\textbf{Proposition 11.5.} If

\[
K(y,x)=\overline{K(x,y)}
\]

almost everywhere, then \(T_\rho\) and \(\widetilde T_\rho\) are
self-adjoint.

\textbf{Proof.} Under the stated symmetry,
\(k_\rho(y,x)=\overline{k_\rho(x,y)}\), so the ambient integral operator
is self-adjoint. The intrinsic result follows from the unitary
intertwining relation. \(\square\)

\section{12. Closure Properties}\label{closure-properties}

\textbf{Proposition 12.1.} Let \(\{\widetilde A_\rho\}\) and
\(\{\widetilde B_\rho\}\) be uniformly bounded, fiber-compatible,
operator-norm Lipschitz families. Then their sum and composition are
also fiber-compatible and operator-norm Lipschitz.

Suppose

\[
\|\widetilde A_\rho\|\le M_A,
\qquad
\|\widetilde B_\rho\|\le M_B,
\]

and their density Lipschitz constants are \(L_A,L_B\). Then the sum has
Lipschitz constant \(L_A+L_B\), while the composition has Lipschitz
constant at most

\[
M_A L_B+M_B L_A.
\]

\textbf{Proof.} Fiber compatibility is preserved under sums and
compositions. For the composition,

\[
\widetilde A_\rho\widetilde B_\rho-
\widetilde A_\sigma\widetilde B_\sigma
=
\widetilde A_\rho(\widetilde B_\rho-\widetilde B_\sigma)
+
(\widetilde A_\rho-\widetilde A_\sigma)\widetilde B_\sigma.
\]

Taking norms yields the stated estimate. \(\square\)

\section{13. Discussion}\label{discussion}

The results separate three levels of structure. First, the intrinsic
operator \(A_\rho:H_\rho\to H_\rho\) determines an observable action
\(A_\rho^\sharp:V_\rho\to V_\rho\). Second, the metric completion
contains the defect subspace \(Z_\rho\), which is absent from the
original weighted fiber. The action on this subspace is not determined
by \(A_\rho\). Third, an extension policy determines how the observable
and defect actions are combined across changing densities.

This distinction explains why the same intrinsic family may induce
different dynamics on the completion. The identity family can produce
the discontinuous support policy \(Q_\rho\) or the constant ambient
identity policy. Continuity is therefore a property of the selected
ambient realization, not of the intrinsic operators alone.

The support projection results show that the \(L^1\)-topology does not
by itself control positivity sets. The support policy is strongly
continuous exactly at densities that are positive almost everywhere, and
it is nowhere continuous in operator norm. Operators depending
discontinuously on support may therefore fail to vary continuously under
\(L^1\)-density convergence. By contrast, the multiplication and
Hilbert-Schmidt constructions show that density-dependent operators can
remain stable when their dependence on the density is mediated through
norm-continuous quantities.

The fixed-density extension classification is elementary at the level of
Hilbert-space block operators. Its role here is to identify precisely
which boundary data remain unspecified after the density-projection
completion. The present results concern bounded operators. A natural
continuation is to study density-dependent closed quadratic forms and
their associated self-adjoint operators, with particular attention to
resolvent and semigroup stability under degeneration of the density.
Such a framework may eventually be applicable to density-dependent
perturbations of Schrödinger-type operators, although domain stability
and lower-semiboundedness would first have to be established. Further
directions include nonlinear maps and evolution equations on the
condensation space. These topics require separate analysis.

\section{14. Conclusion}\label{conclusion}

We studied bounded operator families acting on the varying weighted
Hilbert spaces

\[
H_\rho=L^2(X,\rho\,d\lambda)
\]

and their realizations on the density-projection condensation space

\[
\widehat{\mathcal D}^{DP}\cong L^2(X,\lambda)\times\mathcal R.
\]

The aligned intrinsic operator \(A_\rho^\sharp:V_\rho\to V_\rho\)
determines the action on the observable subspace but not on the
zero-density defect subspace \(Z_\rho\). We showed that the ambient
extension space is affine and modeled on
\(\mathcal B(Z_\rho,\mathcal H)\). The zero-defect extension attains the
minimum possible operator norm, and its structural and spectral
properties were identified.

At the family level, continuity of the induced condensation-space
dynamics was shown to be equivalent to strong-operator continuity of the
selected extension policy. Global Lipschitz continuity was proved to be
rigid: it forces the ambient operator family to be independent of the
density. We also characterized convergence of observable support
projections, identified the exact strong-continuity points of the
support policy, proved that this policy is nowhere operator-norm
continuous, and demonstrated that \(L^1\)-density convergence alone does
not control support-dependent operators.

Finally, smoothed multiplication families and density-weighted
Hilbert-Schmidt operators provided nontrivial stable examples. These
results show that zero-density defects create genuine extension freedom
but do not force instability in every density-dependent model. Stability
depends on how an extension policy incorporates both the magnitude and
the support of the observation density.

\section*{Declaration of Generative AI Assistance}

OpenAI's ChatGPT was used as an assistive tool in the preparation
and revision of this manuscript, including manuscript organization,
language editing, LaTeX preparation, and refinement of the exposition
of mathematical arguments. The author independently reviewed the
final manuscript and takes full responsibility for all mathematical
claims, proofs, references, and conclusions.


\begin{thebibliography}{99}

\bibitem{KuwaeShioya2003}
K.~Kuwae and T.~Shioya,
``Convergence of spectral structures: A functional analytic theory and its applications to spectral geometry,''
\emph{Communications in Analysis and Geometry}, vol.~11, no.~4, pp.~599--673, 2003.
\newblock doi:10.4310/CAG.2003.v11.n4.a1.

\bibitem{MugnoloNittkaPost2013}
D.~Mugnolo, R.~Nittka, and O.~Post,
``Norm convergence of sectorial operators on varying Hilbert spaces,''
\emph{Operators and Matrices}, vol.~7, no.~4, pp.~955--995, 2013.
\newblock doi:10.7153/oam-07-54.

\bibitem{PostZimmer2022}
O.~Post and S.~Zimmer,
``Generalised norm resolvent convergence: Comparison of different concepts,''
\emph{Journal of Spectral Theory}, vol.~12, no.~4, pp.~1459--1506, 2022.
\newblock doi:10.4171/JST/442.

\bibitem{Yetter2005}
D.~N. Yetter,
``Measurable categories,''
\emph{Applied Categorical Structures}, vol.~13, no.~5--6, pp.~469--500, 2005.
\newblock doi:10.1007/s10485-005-9003-6.

\bibitem{DPCompletion}
S.~Jeon,
``Density-Projection Condensation Spaces: A Metric Completion of Function-Density Pairs over Varying Weighted $L^2$-Spaces,''
unpublished manuscript, 2026.

\bibitem{Conway1990}
J.~B. Conway,
\emph{A Course in Functional Analysis}, 2nd ed., Springer, 1990.

\bibitem{Folland1999}
G.~B. Folland,
\emph{Real Analysis: Modern Techniques and Their Applications}, 2nd ed., Wiley, 1999.

\bibitem{ReedSimon1980}
M.~Reed and B.~Simon,
\emph{Methods of Modern Mathematical Physics, Vol.~I: Functional Analysis}, Academic Press, 1980.

\end{thebibliography}
\end{document}